\documentclass[11pt]{amsart}
\usepackage{amsmath}

\theoremstyle{definition}

\numberwithin{equation}{section}
\input{tcilatex}

\begin{document}
\title[Explicit formula]{An Explicit Formula Relating Stieltjes Constants
and Li's Numbers}
\author{Krzysztof Ma\'{s}lanka}
\address{Astronomical Observatory of the Jagiellonian University \\
Orla 171, 30-244 Cracow, Poland}
\date{13 June 2004}
\thanks{e-mail: maslanka@oa.uj.edu.pl}

\begin{abstract}
In this paper we present a new formula relating Stieltjes numbers $\gamma
_{n}$ and Laurent coefficinets $\eta _{n}$ of logarithmic derivative of the
Riemann's zeta function. Using it we derive an explicit formula for the
oscillating part of Li's numbers $\overset{\sim }{\lambda }_{n}$ which are
connected with the Riemann hypothesis.
\end{abstract}

\maketitle

\section{Introduction}

First recall some basic definitions and conventions. The Stieltjes constants 
$\gamma _{n}$ are essentially the coefficients in the Laurent expansion of
the Riemann zeta function about $s=1$:%
\begin{equation}
\zeta \left( s+1\right) =\frac{1}{s}+\sum_{n=0}^{\infty }\gamma _{n}s^{n}
\label{series gamma}
\end{equation}%
(Sometimes another convention is adopted. It differs from the above by the
factor $\left( -1\right) ^{n}/n!$) It is well known that zeta has a single
simple pole at $s=1$ with residue $1$ which is evident from (\ref{series
gamma}).

Other useful coefficients $\eta _{n}$ are these which appear in the Laurent
expansion of the logarithmic derivative of zeta about $s=1$:%
\begin{equation}
-\frac{\zeta ^{\prime }}{\zeta }\left( s+1\right) =\frac{1}{s}%
+\sum_{n=0}^{\infty }\eta _{n}s^{n}  \label{series eta}
\end{equation}%
It can be shown that \cite{BombieriLagarias} 
\begin{equation}
\gamma _{n}:=\frac{(-1)^{n}}{n!}\lim_{x\rightarrow \infty }\left(
\sum_{k\leq x}\frac{(\log k)^{n}}{k}-\frac{(\log x)^{n+1}}{n+1}\right)
\label{gamman}
\end{equation}%
\begin{equation}
\eta _{n}:=\frac{\left( -1\right) ^{n}}{n!}\lim_{x\rightarrow \infty }\left(
\sum_{k\leq x}\Lambda \left( k\right) \frac{(\log k)^{n}}{k}-\frac{(\log
x)^{n+1}}{n+1}\right)  \label{etan}
\end{equation}%
where $\Lambda \left( k\right) $ denotes the von Mangoldt function (see e.g. 
\cite{Edwards}, p. 50). It is related to the Riemann zeta function by%
\begin{equation*}
-\frac{\zeta ^{\prime }\left( s\right) }{\zeta \left( s\right) }%
=\dsum\limits_{n=2}^{\infty }\frac{\Lambda \left( n\right) }{n^{s}}
\end{equation*}%
where $\Re s>1$. It may be shown that $\Lambda \left( n\right) $ is zero
unless $n$ is a prime power: $n=$ $p^{k}$, in which case $\Lambda \left(
n\right) $ is $\log (p)$

However, the expansions (\ref{gamman}) and (\ref{etan}) are slowly
convergent and are therefore useless in numerical computations. Effective
numerical algorithms for calculating $\gamma _{n}$ were given by Keiper \cite%
{Keiper} and\ Kreminski \cite{Kreminski}. An interesting recurrence relation
involving both $\gamma _{n}$ and $\eta _{n}$ has been recently discovered by
Coffey \cite{Coffey} (written here in our convention (\ref{series gamma})
and (\ref{gamman}) concerning $\gamma _{n}$):%
\begin{equation}
\eta _{n}=-\left( n+1\right) \gamma _{n}-\dsum\limits_{k=0}^{n-1}\eta
_{k}\gamma _{n-k-1}  \label{Recurrence}
\end{equation}

\section{Derivation}

Here are several initial relations obtained directly using (\ref{series eta}%
): 
\begin{eqnarray}
\eta _{0} &=&-\gamma _{0}  \label{etas} \\
\eta _{1} &=&+\gamma _{0}^{2}-2\gamma _{1}  \notag \\
\eta _{2} &=&-\gamma _{0}^{3}+3\gamma _{0}\gamma _{1}-3\gamma _{2}  \notag \\
\eta _{3} &=&+\gamma _{0}^{4}-4\gamma _{0}^{2}\gamma _{1}+2\gamma
_{1}^{2}+4\gamma _{0}\gamma _{2}-4\gamma _{3}  \notag \\
\eta _{4} &=&-\gamma _{0}^{5}+5\gamma _{0}^{3}\gamma _{1}-5\gamma _{0}\gamma
_{1}^{2}-5\gamma _{0}^{2}\gamma _{2}+5\gamma _{1}\gamma _{2}+5\gamma
_{0}\gamma _{3}-5\gamma _{4}  \notag \\
&&...  \notag
\end{eqnarray}%
With growing $n$ the number of terms in $\eta _{n}$ increases as the
so-called partition number.

At the first glance there seem to be some regularities in (\ref{etas}).
Indeed, careful (and very tedious) inspection of the indices in $\gamma $s,
their exponents, as well as the numerical coefficients reveals several
simple rules:

\begin{enumerate}
\item Each $\eta _{n}$ is a linear combination of products of powers of
Stieltjes numbers with some numerical coefficient: $\pm A\gamma
_{i}^{k_{i}}\gamma _{j}^{k_{j}}...\gamma _{m}^{k_{m}}$.

\item Every such term appears only once. More precisely, $\pm A\gamma
_{i}^{k_{i}}\gamma _{j}^{k_{j}}...\gamma _{m}^{k_{m}}$ appears in $\eta _{n}$
with $n=(i+1)k_{i}+(j+1)k_{j}+...(m+1)k_{m}$.

\item The sign of each term is determined by the sum of \textit{all}
exponents $k_{i}$: it is equal to $\left( -1\right) ^{p}$, where $%
p=\tsum\limits_{i=0}^{n}k_{i}$.

\item The numerical coefficients $A$ which appear in $\eta _{n}$ may be
obtained as follows. Let us discard for a while in each term power of $%
\gamma _{0}$, i.e. $\gamma _{0}^{k_{0}}$, if there is one. It turns out that 
$A$ contains two factors. The first is%
\begin{equation*}
\frac{n}{k_{i}!k_{j}!...k_{m}!}
\end{equation*}%
and the other is%
\begin{equation*}
\tprod\limits_{i=s+1}^{q-2}\left( n-i\right)
\end{equation*}%
where $s$ is the sum of all indices of $\gamma $ in a given product and $q$
is the number of all $\gamma $s in the same product (in both cases $\gamma $%
s are counted according to their multiplicity and in both cases we neglect $%
\gamma _{0}$), i.e.%
\begin{eqnarray*}
s &=&ik_{i}+jk_{j}+...+mk_{m} \\
q &=&k_{i}+k_{j}+...+k_{m}
\end{eqnarray*}%
The above product may be further simplified:%
\begin{equation*}
\tprod\limits_{i=s+1}^{q-2}\left( n-i\right) =-\left( -1\right) ^{q}\frac{%
\Gamma \left( s+q-n\right) }{\Gamma \left( s-n+1\right) }=\frac{\Gamma
\left( n-s\right) }{\Gamma \left( n-s-q+1\right) }
\end{equation*}%
where the last equality stems from the properties of the Pochhammer's
symbol. According to (2) we have $n-s-q+1=k_{0}+1$, and therefore $\Gamma
\left( n-s-q+1\right) =k_{0}!$ The last step is equivalent to restoring the
powers of $\gamma _{0}$ in the general term $\pm A\gamma _{i}^{k_{i}}\gamma
_{j}^{k_{j}}...\gamma _{m}^{k_{m}}$.
\end{enumerate}

Collecting all the above rules (1-4) together leads to the general
expression for $\eta $ as a function of appropriate $\gamma $s: 
\begin{equation}
\fbox{$\eta _{n-1}=n\dsum\limits_{k_{i}=0}^{n}\left[ \Gamma \left( p\right)
\delta _{n,r}\dprod\limits_{i=0}^{n}\frac{\left( -\gamma _{i}\right) ^{k_{i}}%
}{k_{i}!}\right] $}  \label{EtaFormula}
\end{equation}%
where 
\begin{eqnarray}
p &\equiv &p\left( k_{0},k_{1},...,k_{n}\right) =\sum_{i=0}^{n}k_{i}
\label{Additional} \\
r &\equiv &r\left( k_{0},k_{1},...,k_{n}\right) =\sum_{i=0}^{n}\left(
1+i\right) k_{i}  \notag
\end{eqnarray}%
and the sum (\ref{EtaFormula})\ is performed over \textit{all} combinations
of integers $k_{i}=0,1,2...$ satisfying the constraints (\ref{Additional}).
This sum contains formally many terms, roughly $O(n^{n})$, but, due to the
Kronecker delta which picks out only appropriate terms (i.e. when $r=n$) and
cuts off all the others, the number of non-zero terms is much smaller. This
neat and concise relation (\ref{EtaFormula}) is in fact pretty sophisticated.

Formula (\ref{EtaFormula}) may be inverted. Writing down several initial
expressions (by solving the quasi-linear system (\ref{etas}) with respect to 
$\gamma _{n}$) we get:%
\begin{eqnarray}
\gamma _{0} &=&-\eta _{0}  \label{gammas} \\
\gamma _{1} &=&\frac{1}{2}\left( +\eta _{0}^{2}-\eta _{1}\right)   \notag \\
\gamma _{2} &=&\frac{1}{3!}\left( -\eta _{0}^{3}+3\eta _{0}\eta _{1}-2\eta
_{2}\right)   \notag \\
\gamma _{3} &=&\frac{1}{4!}\left( +\eta _{0}^{4}-6\eta _{0}^{2}\eta
_{1}+3\eta _{1}^{2}+8\eta _{0}\eta _{2}-6\eta _{3}\right)   \notag \\
\gamma _{4} &=&\frac{1}{5!}\left( -\eta _{0}^{5}+10\eta _{0}^{3}\eta
_{1}-15\eta _{0}\eta _{1}^{2}-20\eta _{0}^{2}\eta _{2}+20\eta _{1}\eta
_{2}+30\eta _{0}\eta _{3}-24\eta _{4}\right)   \notag \\
&&...  \notag
\end{eqnarray}%
Again, several regularities are evident. When collected together in the same
way as before they give:%
\begin{equation}
\fbox{$\gamma _{n-1}=\dsum\limits_{k_{i}=0}^{n}\left[ \delta
_{n,r}\dprod\limits_{i=0}^{n}\frac{1}{k_{i}!}\left( \frac{-\eta _{i}}{1+i}%
\right) ^{k_{i}}\right] $}  \label{GammaFormula}
\end{equation}%
where $r$ is given by the second relation (\ref{Additional}). It should be
emphasized that formulas (\ref{GammaFormula}) and (\ref{EtaFormula}), as
well as (\ref{LambdaFormula}) below, may be effectively implemented e.g.
using \textit{Mathematica} symbolic package (it is not trivial since the
number of summations and products is variable):\bigskip 

Table[\{Subscript[\TEXTsymbol{\backslash}[Eta],n-1],

Factor[n*Sum[$\Gamma $[Sum[Subscript[k,i],\{i,0,n\}]]*

KroneckerDelta[n,Sum[(1+i)*Subscript[k,i],\{i,0,n\}]]*

Product[(-Subscript[\TEXTsymbol{\backslash}[Gamma],i])\symbol{94}%
Subscript[k,i]/Subscript[k,i]!,\{i,0,n\}],

Evaluate[Apply[Sequence,Table[\{Subscript[k,j],0,n\},\{j,0,n\}]]]

]]\},\{n,1,5\}]\bigskip

where the Euler gamma function $\Gamma $ above is slightly modified:\bigskip

$\Gamma $[n\_] := Gamma[If[n == 0, 1, n]]

\section{Applications}

It was shown elsewhere by the author \cite{Maslanka}\ that the behavior of
certain numbers $\overset{\sim }{\lambda }_{n}$ with growing $n$ is crucial
for the Riemann hypothesis to be true. Now using%
\begin{equation}
\overset{\sim }{\lambda }_{n}=-\sum_{j=1}^{n}\binom{n}{j}\eta _{j-1}
\label{LambdaOsc}
\end{equation}%
together with (\ref{EtaFormula}) we get%
\begin{equation}
\fbox{$\overset{\sim }{\lambda }_{n}=-\dsum\limits_{k_{i}=0}^{n}\left[
\Gamma \left( p\right) \binom{n}{r}r\dprod\limits_{i=0}^{n}\frac{\left(
-\gamma _{i}\right) ^{k_{i}}}{k_{i}!}\right] $}  \label{LambdaFormula}
\end{equation}%
which explicitly expresses the oscillating part of Li's numbers using
Stieltjes constants. Similarly as in (\ref{EtaFormula}), due to the binomial
coefficient there are in fact many redundant terms in (\ref{LambdaFormula})
which are identically zero, i.e. always when $r$ is greater than $n$. Those
which are non-zero ($r\leq n$) give:

\begin{eqnarray}
\overset{\sim }{\lambda }_{1} &=&\gamma _{0}  \label{lambdas} \\
\overset{\sim }{\lambda }_{2} &=&2\gamma _{0}-\gamma _{0}^{2}+2\gamma _{1} 
\notag \\
\overset{\sim }{\lambda }_{3} &=&3\gamma _{0}-3\gamma _{0}^{2}+\gamma
_{0}^{3}+6\gamma _{1}-3\gamma _{0}\gamma _{1}+3\gamma _{2}  \notag \\
\overset{\sim }{\lambda }_{4} &=&4\gamma _{0}-6\gamma _{0}^{2}+4\gamma
_{0}^{3}-\gamma _{0}^{4}+12\gamma _{1}-12\gamma _{0}\gamma _{1}+4\gamma
_{0}^{2}\gamma _{1}-2\gamma _{1}^{2}  \notag \\
&&+12\gamma _{2}-4\gamma _{0}\gamma _{2}+4\gamma _{3}  \notag \\
&&......  \notag
\end{eqnarray}%
With growing $n$ the number of terms in $\overset{\sim }{\lambda }_{n}$
increases as the summatory function for partition number, i.e. faster than
in the case of (\ref{EtaFormula}).

\section{Discussion}

The distribution of values for various terms which contribute to (\ref%
{lambdas}) is very nonuniform. When $n\rightarrow \infty $ the majority of
values is concentrated near zero. In Figure 1 there are eight histograms for 
$n=3,4,...,10$. The shape of these histograms is roughly symmetrical about
zero but it is just the slight departure from symmetry which contributes to $%
\overset{\sim }{\lambda }_{n}$.

\textit{(Figure 1 about here)}

The fundamental problem, equivalent to the Riemann hypothesis itself is
whether the oscillating part of Li's numbers $\overset{\sim }{\lambda }_{n}$
is bounded by the trend $\overset{-}{\lambda }_{n}$ which grows
asymptotically as:%
\begin{equation*}
\overset{-}{\lambda }_{n}\sim \frac{1}{2}\left( 1+n\ln n\right) +cn
\end{equation*}%
where (cf. \cite{Voros}, \cite{Lagarias}):%
\begin{equation*}
c=\frac{1}{2}\left( \gamma -1-\ln 2\pi \right) 
\end{equation*}

\textbf{Acknowledgment}: I would like to thank dr. Michael Trott, Wolfram
Research, for pointing my attention to a typographical error in a formula as
well as for checking the calculations using \textit{Mathematica}.

\end{document}